\documentclass{birkjour}
\usepackage{amssymb}
\newtheorem{thm}{Theorem}[section]
\newtheorem{cor}[thm]{Corollary}

\theoremstyle{definition}

\theoremstyle{remark}

\numberwithin{equation}{section}

\def\be{\begin{equation}}
\def\ee{\end{equation}}
\newcommand{\bthm}{\begin{thm}}
\newcommand{\ethm}{\end{thm}}
\newcommand{\bcor}{\begin{cor}}
\newcommand{\ecor}{\end{cor}}

\newcommand{\beq}{\begin{eqnarray}}
\newcommand{\beqq}{\begin{eqnarray*}}
\newcommand{\eeq}{\end{eqnarray}}
\newcommand{\eeqq}{\end{eqnarray*}}
\newcommand{\ba}{\begin{array}}
\newcommand{\ea}{\end{array}}
\newcommand{\bee}{\begin{enumerate}}
\newcommand{\eee}{\end{enumerate}}
\newcommand{\bpf}{\begin{proof}}
\newcommand{\epf}{\end{proof}}

\newcommand{\IC}{{\mathbb C}}
\newcommand{\ID}{{\mathbb D}}

\newcommand{\ds}{\displaystyle}

\newcounter{alphabet}

\newenvironment{Thm}[1][]{\refstepcounter{alphabet}%
\bigskip%
\noindent%
{\bf Theorem \Alph{alphabet}}%
{\bf .} \itshape}{\vskip 8pt}

\newenvironment{Lem}[1][]{\refstepcounter{alphabet}%
\bigskip%
\noindent%
{\bf Lemma \Alph{alphabet}}%
{\bf .} \itshape}{\vskip 6pt}

\begin{document}

\title
[Inequality for solutions to nonhomogeneous biharmonic equations]
{On Schwarz-Pick type inequality and Lipschitz continuity
for solutions to nonhomogeneous biharmonic equations}

\author{Peijin Li}
\address{P. Li, Department of Mathematics, Hunan First Normal University, Changsha,
Hunan 410205, People's Republic of China.} \email{wokeyi99@163.com}

\author{Yaxiang Li}
\address{Y. Li, Department of Mathematics, Hunan First Normal University, Changsha,
Hunan 410205, People's Republic of China.}
\email{yaxiangli@163.com}

\author{Qinghong Luo}
\address{Qinghong Luo, Department of Mathematics,
Hunan Normal University, Changsha, Hunan 410081, People's Republic of China
and Department of Mathematics,
Hunan First Normal University, Changsha, Hunan 410205, People's Republic of China}
\email{luoqh207@qq.com}

\author{Saminathan Ponnusamy}
\address{S. Ponnusamy, Department of Mathematics,
Indian Institute of Technology Madras, Chennai-600 036, India, and\\
Lomonosov Moscow State University, Moscow Center of Fundamental and Applied Mathematics, Moscow, Russia.}
\email{samy@iitm.ac.in}

\date{August 06, 2022}

\begin{abstract}
The purpose of this paper is to study the Schwarz-Pick type inequality and the  Lipschitz continuity for the solutions to
the nonhomogeneous biharmonic equation: $\Delta(\Delta f)=g$, where $g:$ $\overline{\ID}\rightarrow\mathbb{C}$ is a continuous
function and $\overline{\ID}$ denotes the closure of the unit disk $\ID$ in the complex plane $\mathbb{C}$. In fact, we
establish the following properties for these solutions: Firstly, we show that the solutions $f$ do not always satisfy the
Schwarz-Pick type inequality
$$\frac{1-|z|^2}{1-|f(z)|^2}\leq C,
$$
where $C$ is a constant. Secondly, we establish a general Schwarz-Pick type inequality of $f$ under certain conditions.
Thirdly, we discuss the Lipschitz continuity of $f$, and as applications, we get the Lipschitz continuity with respect to
the distance ratio metric and the Lipschitz continuity with respect to the hyperbolic metric.
\end{abstract}
\maketitle


\section{Introduction }\label{sec-1}
Let $\mathbb{C}\cong\mathbb{R}^2$ denote the complex plane and for $r>0$, let ${\mathbb D}_r=\{z\in \IC:\,|z|<r\}$. Denote by  $\mathbb {D}:= {\mathbb D} _1$,
the open unit disk, $\mathbb{T}=\partial\mathbb{D}$, the boundary of $\mathbb{D}$, and
$\overline{\ID}=\ID\cup \mathbb{T}$, the closure of $\mathbb{D}$.
Furthermore, we denote by $\mathcal{C}^{m}(\Omega)$ the set of all
complex-valued $m$-times continuously differentiable functions from
$\Omega$ into $\mathbb{C}$, where $\Omega$ stands for a subset of
$\mathbb{C}$ and $m\in\mathbb{N}_0:=\mathbb{N}\cup\{0\}$. In
particular, $\mathcal{C}(\Omega):=\mathcal{C}^{0}(\Omega)$ denotes the
set of all continuous functions in $\Omega$.

Let $\varphi\in {\mathcal C}(\mathbb{T})$, $f^{\ast},g\in\mathcal{C}(\overline{\mathbb{D}})$ and $f\in\mathcal{C}^{4}(\mathbb{D})$.
We consider the following nonhomogeneous biharmonic equations with the Dirichlet boundary values:
\be\label{eq-B1}
\begin{cases} \Delta(\Delta f)=g &
\mbox{ in }\ds  \ID,\\
\ds f_{\overline{z}}=\varphi & \mbox{ on }\ds \mathbb{T},\\
\ds f=f^{\ast} & \mbox{ on }\ds \mathbb{T},
\end{cases}
\ee
where $\Delta f= f_{xx}+f_{yy}  =4f_{z \overline{z}}$
is the {\it Laplacian} of $f$.

In particular, if $g\equiv0$, then the solutions to \eqref{eq-B1} are biharmonic
mappings, see \cite{CPW0, SP}  
and the references therein for certain properties of biharmonic mappings.

The nonhomogeneous biharmonic equation arises in areas of continuum mechanics,
including linear elasticity theory and the solution of Stokes flows (cf. \cite{HS}). 
Chen et al. \cite{CLW} discussed the Schwarz-type lemma, Landau-type theorems and bi-Lipschitz properties for the solutions of
nonhomogeneous biharmonic equations \eqref{eq-B1}. Mohapatra et al. \cite{MWZ} discussed the boundary Schwarz lemma for the solutions of
nonhomogeneous biharmonic equations \eqref{eq-B1}.
The solvability of the biharmonic equations has been studied, for example, in \cite{ATM, MM}.
In this paper we investigate Schwarz-Pick type inequality and Lipschitz continuity for the solutions to the nonhomogeneous biharmonic equations \eqref{eq-B1}.

For $z, w\in\mathbb{D}$, let
$$
G(z,w)=|z-w|^{2}\log\left|\frac{1-z\overline{w}}{z-w}\right|^{2}-(1-|z|^{2})(1-|w|^{2})
$$
and
$$P(z,e^{it})=\frac{1-|z|^{2}}{|1-ze^{-it}|^{2}}
$$
denote the {\it biharmonic Green function} and {\it Poisson kernel}, respectively, where $t\in[0,2\pi]$.

By \cite[Theorem 2]{Be}, we see that all solutions to the biharmonic equations (\ref{eq-B1}) are given by

\beq\label{eq-ch-3}
f(z)&=&\mathcal{P}_{f^{\ast}}(z)+\frac{1}{2\pi}\int_{0}^{2\pi}\overline{z}e^{it}f^{\ast}(e^{it})\frac{1-|z|^{2}}{(1-\overline{z}e^{it})^{2}}dt\\
\nonumber
&&-(1-|z|^{2})\mathcal{P}_{\varphi_{1}}(z)-\frac{1}{8}G[g](z),
\eeq
where $\varphi_{1}(e^{it})=\varphi(e^{it})e^{-it}$,
$$\mathcal{P}_{f^{\ast}}(z)=\frac{1}{2\pi}\int_{0}^{2\pi}P(z,e^{it})f^{\ast}(e^{it})dt, \quad
\mathcal{P}_{\varphi_{1}}(z)=\frac{1}{2\pi}\int_{0}^{2\pi}P(z,e^{it})\varphi_{1}(e^{it})dt
$$
and
$$G[g](z)=\frac{1}{2\pi}\int_{\mathbb{D}}G(z,w)g(w)dA(w).
$$
Here $dA(w)$ denotes the Lebesgue area measure in
$\mathbb{D}$.

\section{Main Results}\label{sec-2}



 Let $\mathcal{B}_H$ denote the set of all complex-valued harmonic mapping from $\mathbb{D}$ into itself. Set
$\mathcal{B}^0_H=\{f\in \mathcal{B}_H:\, f(0)=0\}$.
Pavlovi\'c \cite[Theorem 3.6.1]{Pav1} extended the Heinz' classical lemma in the following general form: if $f\in \mathcal{B}_H $
then one has
\be\label{eq-pav1}
\left|f(z)-\frac{1-|z|^{2}}{1+|z|^{2}}f(0)\right|\leq\frac{4}{\pi}\arctan |z| ~\mbox{ for $z\in\mathbb{D}$}.
\ee
Note that the case $f(0)=0$ is due to Heinz \cite{He}. Moreover, it is easy to see that
$p(t)= \frac{4}{\pi}\arctan t-\frac{2}{\pi}t$
is increasing on $[0,1]$ and thus, $p(t)\leq p(1)$, or equivalently,
\be\label{Lem-C}
\frac{4}{\pi}\arctan t\leq\frac{2}{\pi}(t-1)+1 ~\mbox{ for $t\in [0,1]$}
\ee
(see also \cite[Lemma 2.4]{ZMY}). Now, by using \eqref{eq-pav1} and \eqref{Lem-C}, we get
$$\frac{1-|z|^2}{1-|f(z)|^2}\leq\frac{1-|z|^2}{1-\big(\frac{2}{\pi}(|z|-1)+1\big)^2}=\frac{\pi^2}{4}\frac{1+|z|}{|z|+\pi-1}.$$
Since the function $h(t)=\frac{1+t}{t+\pi-1}$ is increasing for $t\in[0, 1)$,
we can easily derive the following Schwarz-Pick type inequality for $f\in \mathcal{B}^0_H$:
\be\label{ScPick0}
\frac{1-|z|^2}{1-|f(z)|^2}\leq\frac{\pi}{2} ~\mbox{ for $z\in\mathbb{D}$}.
\ee


Let $\mathcal{F}(\mathbb{D}, B)=\{f:\mathbb{D}\rightarrow\mathbb{D}: f(0)=0,|\Delta f| \leq B \cdot|D f|^{2}, f \in \mathcal{C}^{2}\}$,
where $|D f|=|f_z|+|f_{\overline{z}}|$.  Recently in 2015, Kalaj \cite{Ka} proved  that if $f$ is a $K$-quasiconformal mapping and $f\in\mathcal{F}(\mathbb{D}, B)$, then there exists a constant $C(B, K)$ such that
\be\label{ScPick1}
\frac{1-|z|^2}{1-|f(z)|^2}\leq C(B, K) ~\mbox{ for $z\in\mathbb{D}$}.
\ee
In the same article in Kalaj asked whether the quasiconformality assumption can be removed. Recently, Zhong et al. \cite[Theorem 1.3]{ZMY} proved that a mapping in the class  $\mathcal{F}(\mathbb{D}, B)$ does not always enjoy the above Schwarz-Pick type inequality
with $C(B)$ in place of $C(B, K)$ in \eqref{ScPick1}.

In the context of our present study, it is then natural to consider whether the Schwarz-Pick type inequality \eqref{ScPick1} holds with an absolute constant $C$ in place of $C(B, K)$ in \eqref{ScPick1}
for the solutions to \eqref{eq-B1}. Compare with \eqref{ScPick0}.

The first aim of this paper is to give a negative answer to the above question by an example. To be more precise our example conveys the following:

\begin{thm}\label{thm-1}
Suppose that $f\in\mathcal{C}^4(\mathbb{D})\cap \mathcal{B}^0_H$ 
and satisfies \eqref{eq-B1}. Then $f$ does not always enjoy the Schwarz-Pick type inequality
\be\label{S-1}
\frac{1-|z|^2}{1-|f(z)|^2}\leq C ~\mbox{ for $z\in\mathbb{D}$},
\ee
where $C$ is a constant. Furthermore, $f$ does not always satisfy the Poisson differential inequality $|\Delta f| \leq B \cdot|D f|^{2}$.
\end{thm}


Although the answer is negative, one can establish a general Schwarz-Pick type
inequality \eqref{S-1} for the solutions to \eqref{eq-B1} under certain conditions.
For example, in \cite{LP}, the Schwarz-Pick type inequality \eqref{S-1} was obtained for $(K, K')$-quasiconformal self-mappings of $\ID$ satisfying the Poisson
differential inequality $|\Delta f| \leq B \cdot|D f|^{2}$. In \cite{ZMY},  the authors obtained a general Schwarz-Pick type inequality
for the self-mappings of $\mathbb{D}$ satisfying the more general Poisson differential inequality $|\Delta f| \leq a \cdot|D f|^{2}+b$ under certain conditions.
In \cite{KZ},  the authors established that if $f\in \mathcal{B}^0_H$ then  the inequality
$$\frac{1-|z|^q}{1-|f(z)|^q}\leq \frac{\pi}{2} ~\mbox{ for $z\in\mathbb{D}$}
$$
holds for each $q>0$.

The second aim of this paper is to give a general Schwarz-Pick type inequality of the solutions to \eqref{eq-B1} under certain conditions.

\begin{thm}\label{thm-2}
For a given $q\in\{1\}\cup[2,+\infty)$, suppose that $f\in \mathcal{B}_H$, $f(0)=0$ 
and satisfies the biharmonic equations \eqref{eq-B1}, where $f^{\ast}\in\mathcal{C}(\overline{\mathbb{D}})$ is analytic in $\mathbb{D}$.
\bee
\item
For $q=1$, we get
$$\frac{1-|z|}{1-|f(z)|}\leq\frac{1}{\frac{2}{\pi}-4(\|\varphi_1\|_{\infty}+\frac{1}{64}\|g\|_{\infty})};$$
\item
For $q\geq 2$, if
$$(2^{q-1}\cdot q+1)\Big(\|\varphi_1\|_{\infty}+\frac{1}{64}\|g\|_{\infty}\Big)<\frac{2}{\pi},$$
then
$$\frac{1-|z|^q}{1-|f(z)|^q}\leq\frac{1}{\frac{2}{\pi}-(2^{q-1}\cdot q+1)(\|\varphi_1\|_{\infty}+\frac{1}{64}\|g\|_{\infty})},$$
\eee
where $\|\varphi_1\|_{\infty}=\sup_{z\in\mathbb{T}}\{|\varphi_1(z)|\}$ and $\|g\|_{\infty}=\sup_{z\in\mathbb{D}}\{|g(z)|\}$.
\end{thm}

Let $D$ and $\Omega$ be domains in $\mathbb{C}$, and let $L$ be a positive constant. Then a mapping
$f: D\rightarrow \Omega$ is said to be Lipschitz if
$$|f(z_1)-f(z_2)|\leq L|z_1-z_2| ~\mbox{ whenever $z_1, z_2\in\mathbb{D}$}.
$$

In \cite{Pa}, Pavlovi$\acute{c}$ proved that the quasiconformality of harmonic homeomorphisms 
of $\ID$ to itself can be characterized in terms of their bi-Lipschitz continuity. Kalaj \cite{Ka} also proved the Lipschitz continuity of quasiconformal harmonic mappings. Recently,  the Lipschitz continuity of solutions to inhomogeneous biharmonic equation
are established in \cite{CLW, CW}.

The third aim of this paper is to consider the Lipschitz continuity of the solutions to \eqref{eq-B1}.

\begin{thm}\label{thm-3}
Suppose that $f\in \mathcal{B}_H$ and satisfies the biharmonic equations \eqref{eq-B1}, where
 $f^{\ast}\in\mathcal{C}(\overline{\mathbb{D}})$ is an univalent analytic function in $\mathbb{D}$ and $f^{\ast}(e^{it})=e^{i\gamma(t)}$.
Then $f$ is Lipschitz continuous.
\end{thm}

Now, we need some preparations to present our next results.
For a subdomain $G\subset \mathbb{C}$ and for all $z_1$, $z_2\in G$, the distance ratio metric $j_{G}$ is defined as
$$j_{G}(z_1,z_2)=\log\left(1+\frac{|z_1-z_2|}{\min\{\delta_G(z_1),\delta_G(z_2)\}}\right),
$$
where $\delta_G(z)$ denotes the Euclidean distance from $z$ to $\partial G$. The distance ratio metric was introduced by Gehring and Palka \cite{go2} and in the above simplified form by Vuorinen \cite{vo2}.
We say that a mapping $f:\,D\to \Omega$ is Lipschitz continuous with respect to the distance ratio metric if there exists a positive constant $L_1$ such that
$$j_{\Omega}(f(z_1),f(z_2))\leq L_1 j_{D}(z_1,z_2)\ \mbox{for all $z_1,$ $z_2\in D$}.
$$

The hyperbolic distance $d_h(z_1, z_2)$ between the two points $z_1$ and $z_2$ in $\Omega$ is defined by
$$\inf_{\gamma}\Big\{\int_{\gamma}\lambda_{\Omega}(z)|dz|\Big\},
$$
where $\gamma$ runs through all rectifiable curves in $\Omega$ which connect $z_1$ and $z_2$.
It is well known that $\lambda_{\ID}(z)=\frac{2}{1-|z|^2}$.
Similarly, we say that a mapping $f:\,D\to \Omega$ is Lipschitz continuous with respect to the hyperbolic metric if there exists a positive constant $L_2$ such that
$$d_h(f(z_1), f(z_2))\leq L_2d_h(z_1, z_2)\ \mbox{for all $z_1,$ $z_2\in D$}.
$$

Under certain conditions, the subject of a harmonic self-mapping of the unit disk that has
Lipschitz continuity with respect to a given metric has attracted the attention of many researchers.
For example, in \cite{KM}, the authors proved that a $K$-quasiconformal harmonic
mapping from $\ID$ onto itself is bi-Lipschitz with respect to hyperbolic metric. In \cite{LP}, the authors proved that
a $(K, K')$--quasiconformal self-mappings of $\ID$ satisfying the
Poisson differential inequality $|\Delta f| \leq B \cdot|D f|^{2}$ is Lipschitz with respect to the distance ratio metric.

As applications of Theorems \ref{thm-2} and \ref{thm-3}, we give the Lipschitz continuity with respect to the distance ratio metric and Lipschitz continuity with respect to the hyperbolic metric of the solutions to \eqref{eq-B1}, respectively.

\begin{cor}\label{cor-3}
Suppose that $f\in \mathcal{B}^0_H$ and satisfies the biharmonic equations \eqref{eq-B1}, where
$f^{\ast}\in\mathcal{C}(\overline{\mathbb{D}})$ is  an univalent analytic function in $\mathbb{D}$ and $f^{\ast}(e^{it})=e^{i\gamma(t)}$.
\bee
\item
If
$$4\Big(\|\varphi_1\|_{\infty}+\frac{1}{64}\|g\|_{\infty}\Big)<\frac{2}{\pi},
$$
then $f$ is Lipschitz continuous with respect to the distance ratio metric.
%
\item
If
$$5\Big(\|\varphi_1\|_{\infty}+\frac{1}{64}\|g\|_{\infty}\Big)<\frac{2}{\pi},
$$
then $f$ is Lipschitz continuous with respect to the hyperbolic metric.
\eee
\end{cor}

The proofs of Theorems \ref{thm-1}, \ref{thm-2}, \ref{thm-3} and Corollary \ref{cor-3} 
will be presented in Section \ref{sec-3}.

\section{The proofs of main results}\label{sec-3}

\subsection{Proof of Theorem \ref{thm-1}}
Consider $f(z)=2|z|^2z^2-|z|^6z^2$, $z\in\ID$. Clearly, $f(0)=0$.
Since $|f(z)|=|z|^4(2-|z|^4)$ and  the function $h(t)=t(2-t)$ is increasing for $t\in[0, 1)$, it follows easily
that $|f(z)|\leq 1$ in $\ID$. Furthermore, a simple calculation shows that 
\beqq
\begin{cases} \Delta(\Delta f)(z)=-1920z^3\overline{z} &
\mbox{ in }\ds  \ID,\\
\ds f_{\overline{z}}(z)=-z^3& \mbox{ on }\ds \mathbb{T},\\
\ds f(e^{i\theta})=e^{2i\theta} & \mbox{ on }\ds \mathbb{T}.
\end{cases}
\eeqq
 Hence, $f$ is a self-mapping of  $\mathbb{D}$ satisfying the biharmonic equations \eqref{eq-B1}.
With $r=|z|^2$, we find that
$$1-|f(z)|^2=1-|z|^8(2-|z|^4)^2=(1-r^2)^2(1+2r^2-r^4)
$$
and thus,
\beqq \lim_{|z|\rightarrow1^{-}}\frac{1-|z|^2}{1-|f(z)|^2}
&=&\lim_{r\rightarrow1^{-}}\frac{1-r}{(1-r^2)^2(1+2r^2-r^4)} =+\infty.
\eeqq
Therefore, $f$ does not always enjoy the Schwarz-Pick type inequality \eqref{S-1}.

Next, if $0<|z|^4<\frac{2}{3}$, a simple calculation shows that
$$|Df|^2=(|f_z|+|f_{\overline{z}}|)^2=64|z|^6(1 -|z|^4)^2 ~\mbox{ and }~ |\Delta f| = 3\big|2-5|z|^2|
$$
and therefore,
$$\lim_{|z|\rightarrow0^{+}}\frac{|\Delta f|}{|Df|^2}
=+\infty.
$$
That is, $f$ does not satisfy the Poisson differential inequality $|\Delta f| \leq B \cdot|D f|^{2}$. This finishes the proof.
 \qed

\vspace{6pt}

We start with some lemmas which are used in the proof of Theorem \ref{thm-2}.

\begin{Lem}{\rm (\cite[Lemma 2.2]{ZMY})}  \label{Lem-A}
Let $f$ be a harmonic self-mapping of $\ID$ satisfying $|f(0)|<\frac{2}{\pi}$. Then for any
$q\geq1$, the inequality
$$\frac{1-|z|^q}{1-|f(z)|^q}\leq\frac{1}{\frac{2}{\pi}-|f(0)|}$$
holds for every $z\in\ID$.
\end{Lem}

\begin{Lem}{\rm (\cite[Lemma 2.3]{ZMY})}  \label{Lem-B}
For any $0\leq y<1$, $0\leq\varepsilon<1$, $q>1$, we have
$$(y+\varepsilon)^q\leq y^q+2^{q-1}q\varepsilon.$$
\end{Lem}


\subsection{Proof of Theorem \ref{thm-2}}
Since $f^{\ast}\in\mathcal{C}(\overline{\mathbb{D}})$ is analytic in $\mathbb{D}$, we have
$$\frac{1}{2\pi}\int_{0}^{2\pi}\overline{z}e^{it}f^{\ast}(e^{it})\frac{1-|z|^{2}}{(1-\overline{z}e^{it})^{2}}dt
=\frac{1-|z|^{2}}{2\pi}\int_{0}^{2\pi}f^{\ast}(e^{it})\frac{\overline{z}e^{it}}{(1-\overline{z}e^{it})^{2}}dt=0.$$
By \eqref{eq-ch-3}, \eqref{eq-pav1} and $|\mathcal{P}_{\varphi_{1}}|\leq \|\varphi_{1}\|_{\infty}$, we see that
\beqq
|f(z)| &\leq&\left|\mathcal{P}_{f^{\ast}}(z)-\frac{1-|z|^{2}}{1+|z|^{2}}\mathcal{P}_{f^{\ast}}(0)\right|
+(1-|z|^{2})\left|\mathcal{P}_{\varphi_1}(z)-\frac{1-|z|^{2}}{1+|z|^{2}}\mathcal{P}_{\varphi_1}(0)\right|\\
&&+\frac{1-|z|^{2}}{1+|z|^{2}}|\mathcal{P}_{f^{\ast}}(0)-(1-|z|^{2})\mathcal{P}_{\varphi_1}(0)|+\frac{1}{8}|G[g](z)|\\
 &\leq&\frac{4}{\pi}\arctan |z|+(1-|z|^{2})\frac{4}{\pi}\arctan |z|\|\varphi_{1}\|_{\infty}\\
 && +\frac{1-|z|^{2}}{1+|z|^{2}}\Big(|\mathcal{P}_{f^{\ast}}(0)-\mathcal{P}_{\varphi_1}(0)|  +|z|^{2}\|\varphi_{1}\|_{\infty}\Big)+\frac{1}{8}|G[g](z)|.
\eeqq
Using the arguments as in the proof of \cite[Theorem 1.1]{CLW}, we have
\be\label{G-1}
|G[g](z)|\leq \frac{1}{8}\|g\|_{\infty}(1-|z|^{2})^{2}.
\ee
It follows from the assumption $f(0)=0$ and \eqref{G-1} that
$$|\mathcal{P}_{f^{\ast}}(0)-\mathcal{P}_{\varphi_1}(0)|=\Big|\frac{1}{8}G[g](0)\Big|\leq\frac{1}{64}\|g\|_{\infty}.$$
Hence, we have the following estimate
\beq\label{h-1}
|f(z)| &\leq&\frac{4}{\pi} \arctan |z|
+\frac{1-|z|^{2}}{1+|z|^{2}}\left(\frac{1}{64}\|g\|_{\infty}+|z|^{2}\left\|\varphi_{1}\right\|_{\infty}\right) \\ \nonumber
&&+\frac{4}{\pi}\left\|\varphi_{1}\right\|_{\infty}\left(1-|z|^{2}\right) \arctan|z|
+\frac{1}{64}\|g\|_{\infty}\left(1-|z|^{2}\right)^{2}.
\eeq
For the case of $q=1$, by \eqref{h-1} and \eqref{Lem-C}, we have
\beqq
\frac{1-|f(z)|}{1-|z|} &\geq&\frac{2}{\pi}-\frac{1+|z|}{1+|z|^2}\left(\frac{1}{64}\|g\|_{\infty}+|z|^{2}\left\|\varphi_{1}\right\|_{\infty}\right)\\
&&-\frac{4}{\pi}\left\|\varphi_{1}\right\|_{\infty}\left(1+|z|\right)\arctan|z|
-\frac{1}{64}\|g\|_{\infty}(1-|z|^2)(1+|z|)\\
&\geq&\frac{2}{\pi}-4\left(\left\|\varphi_{1}\right\|_{\infty}+\frac{1}{64}\|g\|_{\infty}\right)>0.
\eeqq
For the case of $q\geq 2,$ it follows from the assumption $f(0)=0$, \eqref{G-1} and
$$\left\|\varphi_{1}\right\|_{\infty}+\frac{1}{64}\|g\|_{\infty}
<(2^{q-1}\cdot q+1)\Big(\|\varphi_1\|_{\infty}+\frac{1}{64}\|g\|_{\infty}\Big)<\frac{2}{\pi},
$$
that
$$|\mathcal{P}_{f^{\ast}}(0)|\leq |\mathcal{P}_{\varphi_1}(0)|+\frac{1}{8}|G[g](0)|
\leq \left\|\varphi_{1}\right\|_{\infty}+\frac{1}{64}\|g\|_{\infty}<\frac{2}{\pi}.
$$
Hence, $\mathcal{P}_{f^{\ast}}$ is a harmonic self-mapping of $\mathbb{D}$ satisfying $|\mathcal{P}_{f^{\ast}}(0)|\leq\frac{2}{\pi}$.
By Lemma~A, we have
\be\label{Q-1}
\frac{1-|\mathcal{P}_{f^{\ast}}(z)|^q}{1-|z|^q}\geq \frac{2}{\pi}-|\mathcal{P}_{f^{\ast}}(0)|.
\ee
Then by using \eqref{Q-1}, Lemma~B 
and the estimate
$$|f(z)|\leq|\mathcal{P}_{f^{\ast}}(z)|+(1-|z|^2)\Big(\|\varphi_1\|_{\infty}+\frac{1}{64}\|g\|_{\infty}\Big),$$
we get
\beqq
\frac{1-|f(z)|^q}{1-|z|^q}&\geq& \frac{1-|\mathcal{P}_{f^{\ast}}(z)|^q-2^{q-1}\cdot q(1-|z|^2)\Big(\|\varphi_1\|_{\infty}+\frac{1}{64}\|g\|_{\infty}\Big)}{1-|z|^q}\\
&=& \frac{1-|\mathcal{P}_{f^{\ast}}(z)|^q}{1-|z|^q}-2^{q-1}\cdot q\Big(\|\varphi_1\|_{\infty}+\frac{1}{64}\|g\|_{\infty}\Big)\frac{1-|z|^2}{1-|z|^q}\\
&\geq& \frac{2}{\pi}-|\mathcal{P}_{f^{\ast}}(0)|-2^{q-1}\cdot q\Big(\|\varphi_1\|_{\infty}+\frac{1}{64}\|g\|_{\infty}\Big)\\
&\geq& \frac{2}{\pi}-(2^{q-1}\cdot q+1)\Big(\|\varphi_1\|_{\infty}+\frac{1}{64}\|g\|_{\infty}\Big)>0.
\eeqq
The proof of the theorem is complete.
\qed

\vspace{6pt}
In order to prove Theorem \ref{thm-3}, we need the following results.

\begin{Thm}{\rm (\cite[Theorem 1.2]{Pa})} \label{ThmA}
Let $F(e^{it})=e^{i\gamma(t)}$ be a sense-preserving homeomorphism of $\mathbb{T}$ onto iteslt. If $f=\mathcal{P}_{F}$ is a quasiconformal self-mapping of $\ID$, then it is bi-Lipschitz,
i.e., there is a constant $L<\infty$ such that
$$\frac{1}{L}\leq\left|\frac{f(z_1)-f(z_2)}{z_1-z_2}\right|\leq L\;\;\;\;(z_1, z_2\in\ID),$$
and consequently
$$\frac{1}{L}\leq\frac{1-|f(z)|}{1-|z|}\leq L\;\;\;\;(z\in\ID).$$
\end{Thm}

\begin{Lem}{\rm (\cite[Lemma 2.1]{LW})}  \label{Lem-E}
Let $f$ be a function which is continuously differentiable in $\ID$.
Then $f$ is $L$-Lipschitz continuous with $L>0$ if and only if $|Df(z)|\leq L$ in $\ID$.
\end{Lem}




\subsection{Proof of Theorem \ref{thm-3}}
For $z_1,$ $z_2\in \ID$, we have
$$
|f(z_{2})-f(z_{1})|=\left|\int_{[z_{1},z_{2}]}f_{z}(z)dz+f_{\overline{z}}(z)d\overline{z}\right||\leq |Df(z)||z_1-z_2|,
$$
where $[z_{1},z_{2}]$ stands for the segment in $\ID$ with endpoints $z_{1}$ and $z_{2}$. Hence we only need to estimate $|Df(z)|$.
Since $\mathcal{P}_{f^{\ast}}=f^{\ast}$ is an univalent analytic function from $\overline{\mathbb{D}}$ onto $\overline{\mathbb{D}}$ and $f^{\ast}(e^{it})=e^{i\gamma(t)}$, we know from Theorem~C 
 that, for any $z_1, z_2\in\mathbb{D}$,
$$|\mathcal{P}_{f^{\ast}}(z_1)-\mathcal{P}_{f^{\ast}}(z_2)|\leq L|z_1-z_2|,
$$
where $L$ is a positive constant. By Lemma~D, 
we get that
$$|D\mathcal{P}_{f^{\ast}}(z)|\leq  L.
$$
Since the analyticity of $f^{\ast}$ in $\mathbb{D}$ gives
$$\frac{1}{2 \pi} \int_{0}^{2 \pi} \bar{z}e^{it} f^{\ast}\left(e^{it}\right)\frac{1-|z|^{2}}{\left(1-\bar{z} e^{it}\right)^{2}} dt=0,
$$
we know from \eqref{eq-ch-3} that
\beq\label{eq-p1}
f_{z}(z)&=&[\mathcal{P}_{f}(z)]_{z}
+\overline{z}\mathcal{P}_{\varphi_{1}}(z)
-(1-|z|^{2})[\mathcal{P}_{\varphi_{1}}(z)]_{z}\\ \nonumber
&&-\frac{1}{16\pi}\int_{\mathbb{D}}g(w)G_{z}(z,w)dA(w)
\eeq
and
\beq\label{eq-p2}
f_{\overline{z}}(z)&=&[\mathcal{P}_{f}(z)]_{\overline{z}}
+z\mathcal{P}_{\varphi_{1}}(z)
-(1-|z|^{2})[\mathcal{P}_{\varphi_{1}}(z)]_{\overline{z}}\\ \nonumber
&&-\frac{1}{16\pi}\int_{\mathbb{D}}g(w)G_{\overline{z}}(z,w)dA(w).
\eeq
By using \cite[Lemma 2.5]{LP2}, we get
\be\label{eq-p3}
\frac{\|g\|_{\infty}}{16 \pi} \int_{\mathbb{D}}\left(\left|G_{z}(z, w)\right|+\left|G_{\bar{z}}(z, w)\right|\right) d A(w) \leq \frac{23}{48}\|g\|_{\infty}.
\ee
And by the Schwarz-Pick type lemma for harmonic mappings (see \cite{Co}), we have that
\be\label{eq-p4}
|D\mathcal{P}_{\varphi_{1}}(z)| \leq \frac{4}{\pi} \frac{\|\varphi_1\|_{\infty}}{1-|z|^{2}}.
\ee
Then we conclude from \eqref{eq-p1}, \eqref{eq-p2}, \eqref{eq-p3} and \eqref{eq-p4} that
\be\label{eq-p5}
|Df(z)|\leq |D\mathcal{P}_{f^{\ast}}(z)|+2|z||\mathcal{P}_{\varphi_{1}}(z)|+(1-|z|^{2})|D\mathcal{P}_{\varphi_{1}}(z)|+\frac{23}{48}\|g\|_{\infty}\leq M,
\ee
where
$$
M=L+\big (2+\frac{4}{\pi}\big )\|\varphi_1\|_{\infty}+\frac{23}{48}\|g\|_{\infty}.
$$
The proof of the theorem is complete.
\qed

\subsection{Proof of Corollary \ref{cor-3}}
(1) From the hypotheses of Corollary \ref{cor-3}(1) and Theorem \ref{thm-2}, we obtain that
$$\frac{\delta_{\ID}(z)}{\delta_{\ID}(f(z))}=\frac{1-|z|}{1-|f(z)|}\leq \frac{1}{\frac{2}{\pi}-4(\|\varphi_1\|_{\infty}+\frac{1}{64}\|g\|_{\infty})}.
$$
Moreover, from \eqref{eq-p5}, we see that there exists a constant $M$ such that $|Df|\leq M$. Now, we choose an appropriate constant $M_1$
satisfying $M_1= \max\{M, \frac{2}{\pi}\}$ such that $|Df|\leq M_1$. And then, by Lemma~D, 
we have
$$|f(z_{1})-f(z_{2})|\leq M_1|z_1-z_2|.$$
Consequently,
using the Bernoulli inequality, for any two points $z_1$ and $z_2$ in $\mathbb{D}$, we get
\beqq
j_{\ID}(f(z_1), f(z_2))&=&\log\left(1+\frac{|f(z_1)-f(z_2)|}{\min\{\delta_{\ID}(f(z_1)),\delta_{\ID}(f(z_2))\}}\right)\\
&\leq&\log\left(1+\frac{M_1}{\frac{2}{\pi}-4(\|\varphi_1\|_{\infty}+\frac{1}{64}\|g\|_{\infty})}\cdot
\frac{|z_1-z_2|}{\min\{\delta_{\ID}(z_1),\delta_{\ID}(z_2)\}}\right)\\
&\leq&\frac{M_1}{\frac{2}{\pi}-4\left(\|\varphi_1\|_{\infty}+\frac{1}{64}\|g\|_{\infty}\right)}j_{\ID}(z_1, z_2).
\eeqq
The proof of the first part of the corollary is complete.\\

%
(2) In this case, by  
Theorem \ref{thm-2}, we obtain that
$$\frac{1-|z|^2}{1-|f(z)|^2}\leq \frac{1}{\frac{2}{\pi}-5(\|\varphi_1\|_{\infty}+\frac{1}{64}\|g\|_{\infty})}.
$$
This inequality together with \eqref{eq-p5} imply that
$$\frac{|Df(z)|(1-|z|^2)}{1-|f(z)|^2}\leq \frac{M}{\frac{2}{\pi}-5(\|\varphi_1\|_{\infty}+\frac{1}{64}\|g\|_{\infty})}:=M'.
$$

Let $\gamma$ be the hyperbolic geodesic connecting $z_1$ and $z_2$. Then for any $z\in \gamma$,
we obtain that
$$d_h(f(z_1), f(z_2))\leq\int_{f(\gamma)}\lambda_{\mathbb{D}}(w)|dw|\leq\int_{\gamma}(\lambda_{\mathbb{D}}\circ f)(z)|Df(z)||dz|
\leq  M'd_h(z_1, z_2),
$$
where $w=f(z)$. This completes the proof of the second part of the corollary.
\qed

\vspace{8pt}

%

\subsection*{Acknowledgments} The research was partly supported by the Natural Science Foundation of Hunan Province of China (No. 2022JJ40112),
Scientific Fund of Hunan Provincial Education Department (No. 20B118) and NSF of Hebei Science (No. 2021201006).

%


\subsection*{Conflict of Interests}
The authors declare that there is no conflict of interests regarding the publication of this paper.

\subsection*{Data Availability Statement}
The authors declare that this research is purely theoretical and does not associate with any datas.

\normalsize

\end{document}